\documentclass[12pt]{article}
\usepackage{amssymb,amsfonts,amsmath,amsthm, breqn, color, float, mathtools, booktabs, hvfloat, url, mathrsfs,cite}
\usepackage{adjustbox}
\usepackage{caption}
\newcommand{\cc}{{\mathbb C}}

\newcommand{\calj}{{\mathcal J}}

\newcommand{\calr}{{\mathcal R}}

\newcommand{\calf}{{\mathcal F}}

\newcommand{\beq}{\begin{eqnarray*}}
    \newcommand{\feq}{\end{eqnarray*}}
\newcommand{\beqn}{\begin{eqnarray}}
\newcommand{\feqn}{\end{eqnarray}}

\newtheorem{theorem}{Theorem}
\newtheorem*{conj*}{Conjecture}
\makeatletter \@addtoreset{theorem}{section}\makeatother

\newcommand{\nn}{{\mathbb N}}
\makeatletter \@addtoreset{theorem}{section}\makeatother

\makeatletter \@addtoreset{theorem}{section}\makeatother

\newtheorem*{theorem*}{Theorem}

\newtheorem{corollary}[theorem]{Corollary}

\usepackage{algorithm}
\usepackage{graphicx}
\usepackage[noend]{algpseudocode}
\makeatletter
\def\BState{\State\hskip-\ALG@thistlm}
\makeatother

\DeclareCaptionLabelFormat{noname}{#2}
\newlength\myindent
\title{Fixed points of a random restricted growth sequence}

\author{Toufik~Mansour\thanks{ Department of Mathematics, University of Haifa, 199 Abba Khoushy Ave, 3498838 Haifa, Israel;
        \newline e-mail: tmansour@univ.haifa.ac.il}
    \and
    Reza~Rastegar\thanks{Occidental Petroleum Corporation, Houston, TX 77046 and Departments of Mathematics and Engineering, University of Tulsa, OK 74104, USA - Adjunct Professor; e-mail:  reza\_rastegar2@oxy.com}}

\begin{document}
\maketitle
\begin{abstract}
We call $i$ a  fixed point of a given sequence if the value of that sequence at the $i$-th position coincides with $i$. Here, we enumerate fixed points in the class of restricted growth sequences. The counting process is conducted by calculation of generating functions and leveraging a probabilistic sampling method.
\end{abstract}


\section{Introduction}

For any given sequence $\pi$ of length $n$, $i\in [n]:=\{1,\cdots, n\}$ is a fixed point of $\pi$ if the $i$-th entry $\pi_i$ of the sequence $\pi$ is equal to  $i$. We denote the number of fixed points of the sequence $\pi$ by $\calf_\pi$. The term {\it fixed point} is motivated naturally by the same concept in the class of permutations, where it represents any point that is not moved by a permutation. Fixed points and derangements of permutations are well studied due to their importance in various branches of mathematics including algebra, probability, and combinatorics; see for instance \cite{AT, Bo, B,C,DE, HX, DFG, S,W} for a few examples. Recently, a new line of research toward extending results with regards to fixed points in other classes of discrete sequences has emerged; for instance, Archibald, Blecher, and Knopfmacher \cite{ABK} considered fixed points in compositions and words over the alphabet $[k]$. Inspired by them, in this note, we further investigate fixed points for another important class of sequences in combinatorics; namely, restricted growth sequences. These sequences are of interest in connection with set partitions \cite{Mb}, $q$-analogues \cite{maybe1}, certain combinatorial matrices \cite{galvin}, and Gray codes \cite{conflitti1}.

Before we state our results, a few definitions are in order. Throughout this note, we use $\nn$ as the set of all natural numbers. A sequence of natural numbers $\pi=\pi_1\pi_2\cdots \pi_n\in \nn^n$ is called a restricted growth sequence if
\beq
\pi_1 = 1\qquad \mbox{\rm and}  \qquad \pi_{j+1} \leq 1 + \max\{\pi_1,\cdots,\pi_j\} \ \ \ \text{for all} \ 1\leq j <n.
\feq
There is a bijective connection between these sequences and canonical set partitions.
A \textit{partition} of a set $A$ is a collection of non-empty, mutually disjoint subsets, called \textit{blocks}, whose union is the set $A$. A partition $\Pi$ with $k$ blocks is called a \textit{$k$-partition} and denoted by $\Pi = A_1|A_2|\cdots|A_k$. A $k$-partition $A_1|A_2|\cdots |A_k$ is said to be in the \textit{standard form} if the blocks $A_i$ are labeled in such a way that
\beq
\min A_1 < \min A_2<\cdots< \min A_k.
\feq
The partition can be represented equivalently by the \textit{canonical sequential form} $ \pi_1\pi_2\ldots\pi_n,$ where
$\pi_i\in [n]$ and $i\in A_{\pi_i}$ for all $i$ \cite{Mb}. In words, $\pi_i$ is the label of the partition block that contains $i.$
It is easy to verify that a word $\pi\in [k]^n$ is a canonical representation of a $k$-partition of $[n]$ in the standard form
if and only if it is a restricted growth sequence \cite{Mb}.

Throughout this note, we use the terms {\it restricted growth sequence} and {\it set partition} interchangeably. As it becomes clear, in the context of restricted growth sequences, the fixed points are also closely related to the records. We remind the reader that the $i$-th entry $\pi_i$ in the sequence $\pi$ is a record if $\pi_i > \pi_j$
for all $j \in [i-1]$.  Clearly, for any given fixed point $i$ in a restricted growth sequence $\pi$, $i$ is a record. In addition, each $j\in [i-1]$ is a fixed point and hence a record. Therefore, $\calf_\pi$ is precisely the length of the maximal prefix of $\pi$ whose elements are records. We refer the reader to \cite{CMS, KMW, RAT} for a few discussions around records for restricted growth sequences.

We denote by $\calr_n$ the set of all restricted growth sequences of length $n$, and denote by $\calr_{n,k}$ the set of all restricted growth sequences of length $n$ with maximal letter $k$.
\par
\begin{figure}[h!b]
	\centering
	\begin{minipage}{.5\textwidth}
		\centering
		\includegraphics[scale=.4]{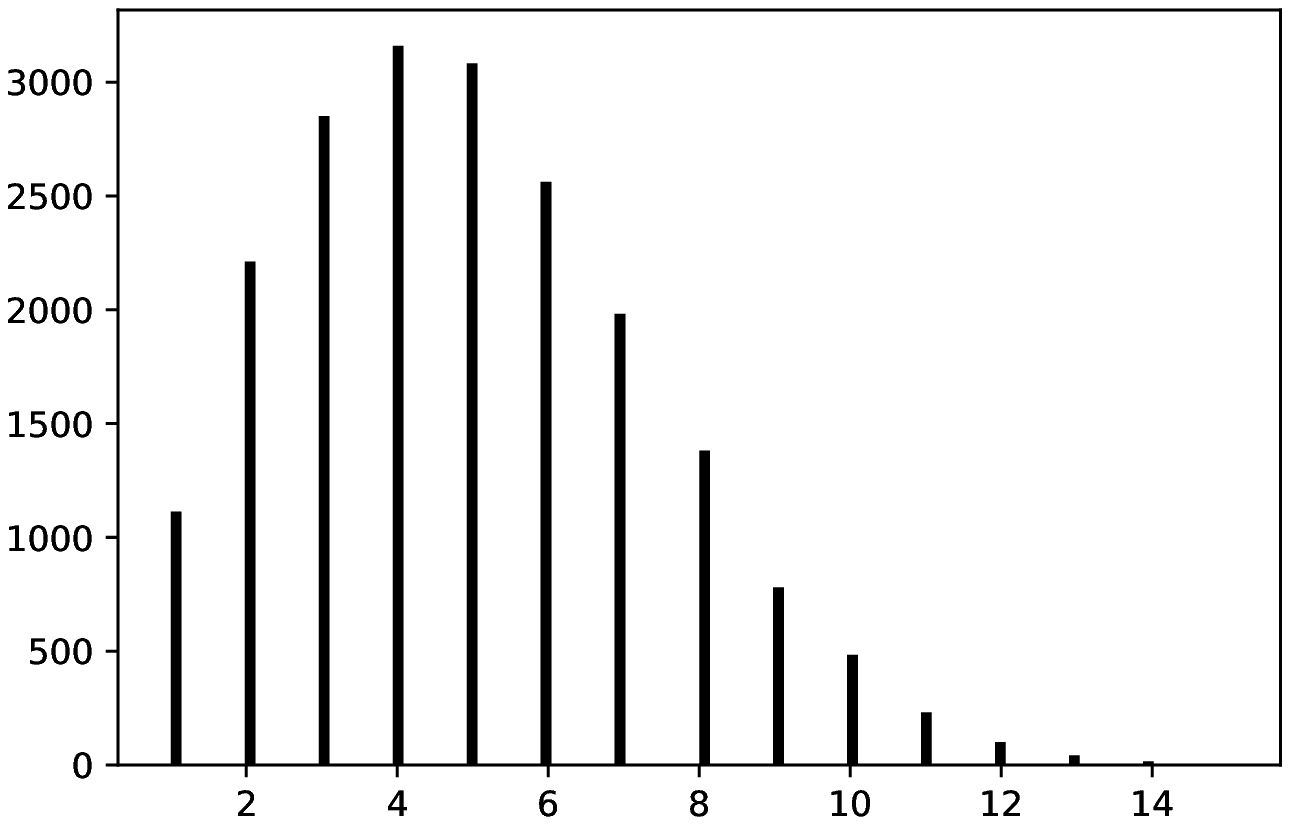}
	\end{minipage}%
	\begin{minipage}{.5\textwidth}
		\centering
		\includegraphics[scale=0.4]{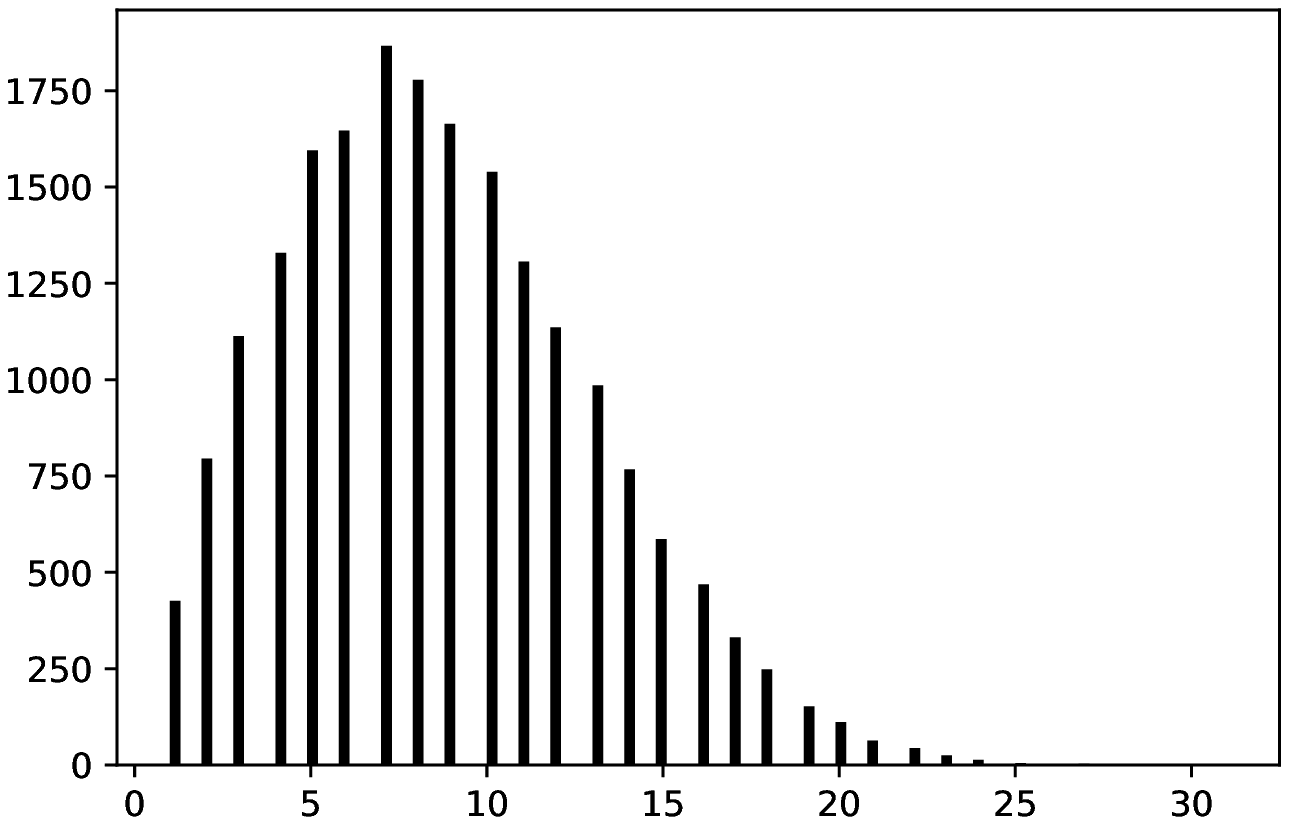}
	\end{minipage}
	\caption{Empirical distributions of $\calf_{50}$ (left) and  $\calf_{200}$ (right) based on $20000$ samples. }
	\label{fig:hn_dist}
\end{figure}

Let $S_{n,k}$ be a Stirling number of the second kind and $B_n$ be the $n$-th Bell number \cite{Mb}. It is well-known that the cardinality of the set $\calr_n$ is $B_n$. In addition, the cardinality of the set $\calr_{n,k}$ is $S_{n,k}$ with the exponential generating function $e^{y(e^x-1)}$, where $y$ counts the number of blocks. The sequence of Bell numbers $(B_n)_{n\geq 0}$ can be then defined, for instance, through the formula
$B_n=\sum_{k=0}^n S_{n,k},$ or, recursively via the formula $B_{n+1}=\sum_{k=0}^n \binom{n}{k}B_k$
with $B_0=1,$ or through Dobinski's formula \cite{comtet}
\beqn
\label{doob}
B_n=\frac{1}{e}\sum_{m=0}^{\infty}\frac{m^n}{m!}, \qquad\qquad n\geq 0.
\feqn
In what follows, we denote a random restricted growth sequence, sampled uniformly from $\calr_{n,k}$ (resp. $\calr_{n}$) by $\pi_{(n,k)}$ (resp. $\pi_{(n)}$). That is,
\beq
P(\pi_{(n,k)} = \pi) = \frac{1}{S_{n,k}} \quad \mbox{for all} \quad \pi \in \calr_{n,k},
\feq
and
\beq
P(\pi_{(n)}  = \pi) = \frac{1}{B_n} \quad \mbox{for all} \quad \pi \in \calr_{n}.
\feq
We use $E(\cdot)$ to refer to the expectation with respect to the probability distribution $P(\cdot)$.
We denote by $\calf_n:=\calf_{\pi_{(n)}}$ the number of fixed points in a uniformly sampled random restricted growth sequence $\pi_{(n)}$. Fig.~\ref{fig:hn_dist} shows  the empirical distributions of $\calf_n$ over $20000$ independently sampled instances of $\pi_{(n)}$ for $n=50$ and $n=200$.

\subsection{Statement of results}
For any $k\in\nn,$ we define $Q_k(x;q)$ to be the ordinary generating function enumerating fixed points over the set of restricted growth sequences in $\cup_{n=k}^\infty \calr_{n,k}$; that is,
\beqn
Q_k(x,q) := \sum_{n=k}^\infty x^n S_{n,k} E\big(q^{\calf_n}|\pi_{(n)}\in\calr_{n,k}\big)
= \sum_{n=k}^\infty  \sum_{\pi \in \calr_{n,k}} x^n q^{\calf_\pi}, \quad x, q\in \cc, \label{Pkxq-def}
\feqn
where $\cc$ is the set of all complex numbers. Knowing an explicit form of \eqref{Pkxq-def}, would in principle give us the distribution of $\calf_n$ in full details for all $n\in\nn.$ However, finding the explicit form of $Q_k(x,q)$ is a daunting task. Hence, we instead study the following exponential generating function:
\beqn \label{R_def}
R(x,y;q):=\sum_{n\geq0}([x^n]Q(x,y;q))\frac{x^n}{n!} = \sum_{k\geq1} y^k \sum_{n\geq 0} \frac{x^n}{n!} B_n E(\calf_n|\pi_{(n)}\in\calr_{n,k}),
\feqn
where $[x^n]f$ denotes the coefficient of $x^n$ in $f$ and $Q(x,y;q)$ is the ordinary generating function for $Q_k(x;q)$; that is,
$$Q(x,y;q):=\sum_{k\geq0}Q_k(x;q)y^k.$$
Our first result states that
\begin{theorem}\label{tha1}
	The exponential generating function $R(x,y;q)$ is given by
	$$R(x,y;q)=e^{y(e^x-1)}+\int_0^x(q-1)ye^{x-t+(1+qt)ye^{x-t}-y}dt.$$
\end{theorem}

The proof of this theorem is given in Section~\ref{proof1}. Note, by Theorem \ref{tha1}, $R(x,y;q)$ reduces to $e^{y(e^x-1)}$ when $q=1$, which is exactly the exponential generating function for the size of $\calr_{n,k}$ (see \cite{Mb}). Moreover, for $q=0$, it implies $$R(x,y;0)=e^{y(e^x-1)}-\int_0^xye^{y(e^t-1)+t}dt=e^{y(e^x-1)}-e^{y(e^x-1)}+1=1,$$
which coincides with the fact that the only restricted growth sequence with no fixed point is the null sequence.

Recall that $\calr_{n,0}$ is an empty set for $n\geq 0$. Let $R_m(x,y)$ be the exponential generating function for the number of sequences in $\calr_{n,k}$ with exactly $m$ fixed points; that is,
\beq
R_m(x,y) :=\sum_{n\geq0}\sum_{k=0}^n \sum_{\substack{\pi\in\calr_{n,k}\\ \calf_\pi=m}}x^ny^k,
\feq
with $R_0(x,y)=1$. Then, Theorem \ref{tha1} implies that
\begin{corollary}
	For all $m\geq1$,
	\begin{align*}
	R_m(x,y)&=\frac{1}{(m-1)!}\int_0^xy^mt^{m-1}e^{y(e^{x-t}-1)+m(x-t)}dt
	-\frac{1}{m!}\int_0^xy^{m+1}t^me^{y(e^{x-t}-1)+(m+1)(x-t)}dt.
	\end{align*}
\end{corollary}

Next, we study the exponential generating function $T(x)$ for the total number of fixed points in $\calr_n$; that is
\beqn \label{Tn}
T(x):=\frac{\partial}{\partial q}R(x,1;q)\mid_{q=1},
\feqn
where by an inductive argument, provided in Section~\ref{proof1}, we conclude that

\begin{corollary} \label{thm2}
	The average number of fixed points over all the restricted growth sequences is
	\beqn \label{Ef_n_gen}
	E(\calf_n) = \frac{1}{B_n}\frac{d^n}{dx^n}T(x)\mid_{x=0}
	=\frac{1}{B_n}\sum_{i=1}^n\left(\sum_{j=1}^i\sum_{\ell=j}^i\frac{(-1)^{\ell-i}\ell^{n-j}}{(\ell-j)!(i-\ell)!}\right).
	\feqn
\end{corollary}

Our final result, Theorem~\ref{thm3} provides a closed form expression for the probability distribution of $\calf_n$ in terms of explicit polynomials
of Bell numbers. We will use the following well-known extension of Dobinski's identity \eqref{doob} to express the result.
Recall that, for any integers $n,t\geq 0$ we have:
\beqn
\nonumber
\Theta_n(t)&:=&\frac{1}{e}\sum_{m=t}^{\infty} \frac{m^n}{(m-t)!} =
\frac{1}{e}\sum_{k=0}^{\infty} \frac{(k+t)^n}{k!} =
\frac{1}{e}\sum_{\ell=0}^n \binom{n}{\ell} t^{n-\ell}\sum_{k=0}^{\infty} \frac{k^\ell}{k!}
\\
&=&
\sum_{\ell=0}^n \binom{n}{\ell} t^{n-\ell} B_\ell, \label{Blextension}
\feqn
where for the last equality we applied the original formula \eqref{doob}. The theorem states
\begin{theorem} \label{thm3}
For $n\in \nn$, we have
\beq
P(\calf_n=j) = \left\{
\begin{array}{lr}
	\frac{B_{n-1}}{B_n} & \text{if } j=1\\
	\frac{j\Theta_{n-j-1}(j)}{B_n} & \qquad \text{if } 1<j<n \\
	\frac{1}{B_n} & \text{if } j=n
\end{array}
\right.
\feq
In particular, for $n\geq 2$, we have
\beqn \label{Ef_n_prob}
E(\calf_n) = \frac{n+\sum_{j=1}^{n-1}j^2\Theta_{n-j-1}(j)}{B_n}.
\feqn
\end{theorem}
\par
The proof of Theorem~\ref{thm3} is given in Section~\ref{pr2}. It is based on the sampling method devised in \cite{stam1} and has been recently exploited for the enumeration of other complex quantities over $\calr_n$. See for instance \cite{RAT2} where the authors enumerate horizontal visibility graphs of sequences in $\calr_n$. We remark that the $j=1$ and $j=n$ cases are obvious; if $\pi\in \calr_n$ such that $\calf_\pi=1$, then it must be that $\pi_2=1$ and clearly there are $B_{n-1}$ choices for such words. Also, there is exactly one word $\pi=12\cdots n\in \calr_n$ for which $\calf_\pi=n$. However, we include our brute-force probabilistic calculation in Section~\ref{pr2} for the sake of completeness.\\

We make a final remark that from \eqref{Ef_n_gen} and \eqref{Ef_n_prob}, we obtain the following identity for $n\geq 2$  which may be of interest independently to the reader:
\beq
\sum_{j=1}^{n-1} \sum_{\ell=0}^{n-j-1} \binom{n-j-1}{\ell} j^{n-j+1-\ell} B_\ell
=\sum_{i=1}^n\left(\sum_{j=1}^i\sum_{\ell=j}^i\frac{(-1)^{\ell-i}\ell^{n-j}}{(\ell-j)!(i-\ell)!}\right) - n.
\feq

\section{Proof of Theorems~\ref{tha1} and Corollary~\ref{thm2}}
\label{proof1}

The proof is based on the observation that each restricted growth sequence $\pi$ with maximal value $k$ and $\calf_\pi=i$ fixed points falls into one of the following cases:
\begin{itemize}
	\item[] {\bf Case $\calf_\pi=i=k$}: Here $\pi=12\cdots k\pi^{(k)}$, where $\pi^{(k)}$ is a word over alphabet $[k]$.
	\item[] {\bf Case $\calf_\pi=i<k$}: Here $\pi=12\cdots i\pi^{(i)}(i+1)\pi^{(i+1)}\cdots k\pi^{(k)}$, where $\pi^{(j)}$ is a word over alphabet $[j]$, for all $j=i,i+1,\ldots,k$, and $\pi^{(i)}$ is not the empty word.
\end{itemize}
By considering these cases we arrive at the following equation:
$$Q_k(x;q)=\frac{x^kq^k}{1-kx}
+\sum_{i=1}^{k-1}\frac{ix^{k+1}q^i}{\prod_{j=i}^k(1-jx)},$$
which implies
\begin{align}\label{eqA1}
(1-kx)Q_k(x;q)-xQ_{k-1}(x;q)=x^kq^{k-1}(q-1),
\end{align}
with $Q_0(x;q)=1$ and $Q_1(x;q)=\frac{xq}{1-x}$.

By multiplying \eqref{eqA1} by $y^k$ and summing over $k\geq2$, we obtain
\begin{align*}
&Q(x,y;q)-1-\frac{xyq}{1-x}\\
&-xy\frac{\partial}{\partial y}\left(Q(x,y;q)-1-\frac{xyq}{1-x}\right)-xy(Q(x,y;q)-1)=\sum_{k\geq2}x^kq^{k-1}(q-1)y^k,
\end{align*}
which leads to
\begin{align}\label{eqA2}
(1-xy)Q(x,y;q)-xy\frac{\partial}{\partial y}Q(x,y;q)=\frac{1-xy}{1-xyq}.
\end{align}
Next, we translate \eqref{eqA2} in terms of exponential generating function $R(x,y;q)$ defined by \eqref{R_def} and write
\beq
\frac{\partial}{\partial x}R(x,y;q)-yR(x,y;q)-y\frac{\partial}{\partial y}R(x,y;q)=y(q-1)e^{xyq}.
\feq
Solving this equation with the initial condition $R(0,y;q)=R(x,0;q)=1$ completes the proof of Theorem \ref{tha1}.
\medskip

Recall \eqref{Tn}. Theorem \ref{tha1} shows that $T(x)=\int_0^xe^{(x-t)e^t+e^t+t-1}dt$.
By induction on $m$, we have
\begin{align}\label{eqT1}
\frac{d^m}{dx^m}T(x)=\sum_{i=1}^ma_{mi}e^{e^x+ix-1}+\int_0^xe^{(x-t)e^t+e^t+(m+1)t-1}dt,
\end{align}
where $a_{mi}$ are natural numbers satisfying the relations
\begin{align}\label{eqR1}
\left\{\begin{array}{l}
a_{mm}=1+a_{(m-1)(m-1)},\\
a_{mi}=a_{(m-1)(i-1)}+ia_{(m-1)i},\quad i=1,2,\ldots,m-1,
\end{array}\right.
\end{align}
with $a_{11}=1$. Let $a_i(x)=\sum_{m\geq i}a_{mi}x^m$. Then, by \eqref{eqR1}, we have
\begin{align*}
a_i(x)&=\sum_{j=1}^i\frac{x^i}{(1-jx)(1-(j+1)x)\cdots(1-ix)}.
\end{align*}
Note that, by a partial fraction decomposition, we have
\begin{align*}
\frac{1}{(1-jx)(1-(j+1)x)\cdots(1-ix)}&=\sum_{\ell=j}^i\frac{(-1)^{\ell-i}\binom{i-j}{\ell-j}\ell^{i-j}}{(i-j)!(1-\ell x)}.
\end{align*}
Hence,
\begin{align*}
a_i(x)&=\sum_{j=1}^i\sum_{\ell=j}^i\frac{(-1)^{\ell-i}\binom{i-j}{\ell-j}\ell^{i-j}x^i}{(i-j)!(1-\ell x)}.
\end{align*}
We inspect the coefficient of $x^m$ in $a_i(x)$ and derive
$$a_{mi}=\sum_{j=1}^i\sum_{\ell=j}^i\frac{(-1)^{\ell-i}\binom{i-j}{\ell-j}\ell^{m-j}}{(i-j)!}.$$
Therefore, \eqref{eqT1} yields
\begin{align*}
\frac{d^m}{dx^m}T(x)=\sum_{i=1}^m\left(\sum_{j=1}^i\sum_{\ell=j}^i\frac{(-1)^{\ell-i}\binom{i-j}{\ell-j}\ell^{m-j}}{(i-j)!}\right)e^{e^x+ix-1}+\int_0^xe^{(x-t)e^t+e^t+(m+1)t-1}dt.
\end{align*}
Finally, setting $x=0$, we complete the proof of Corollary \ref{thm2}.

\section{Proof of Theorem~\ref{thm3}}
\label{pr2}
The proof relies on the use of a generator of a uniformly random set partition of $[n]$ proposed by Stam \cite{stam1}.
We next describe Stam's algorithm for a given $n.$
\begin{enumerate}
\item For $m\in\nn,$ let $\mu_n(m)=\frac{m^n}{m!\ e B_n}.$ Dobinski's formula \eqref{doob} shows that $\mu_n(\,\cdot\,)$ is a probability distribution on $\nn$.
\par
At time zero, choose a random $M\in\nn$ distributed according to $\mu_n,$ and arrange $M$ empty and unlabeled boxes.
\item Arranges $n$ balls labeled by integers from the set $[n].$
\par
At time $i\in [n],$ place the ball `$i$' into one of the $M$ boxes, chosen uniformly at random. Repeat until there are no balls remaining.
\item Label the boxes in the order that they cease to be empty. Once a box is labeled, the label does not change.
\item Form a set partition $\pi$ of $[n]$ with $i$ in the $k$-th block if and only if ball `$i$" is in the $k$-th box.
\end{enumerate}
Let $N_i$ be the random number of nonempty boxes right after placing the $i$-th ball and $X_i$ be the label
of the box where the $i$-th ball was placed. Notice that if the $i$-th ball is dropped in an empty box, then
$X_i=N_{i-1}+1$ and $N_{i} = N_{i-1}+1.$ Otherwise, if the box was occupied previously, $X_i=X_j$ where $j<i$ is the first ball that was dropped in that box and $N_i = N_{i-1}$. Then, $X:=X_1\cdots X_n$ is the random set partition of $[n]$ produced  by the algorithm.
\par
We denote by $P_m(\,\cdot\,)$ the conditional probability distribution $P(\,\cdot\,|\,M=m).$ Clearly $N_1=1$, $N_i\leq i,$ and
\beq
P_m(N_{i+1}=t+1|N_i=t) = \frac{m-t}{m}
\qquad
\mbox{\rm and}
\qquad
P_m(N_{i+1}=t|N_i=t) = \frac{t}{m}.
\feq
Let $\alpha_{i,t}(m):=P_m(N_i=t).$ Then, taking into account that
\beq
P_m(N_i=t) = P_m(N_i=t, N_{i-1}=t-1) + P_m(N_i=t, N_{i-1}=t),
\feq
we obtain:
\beq
\alpha_{i,t} (m)=
\left\{
\begin{array}{ll}
\frac{t}{m} \alpha_{i-1,t}(m) + \frac{m-t+1}{m}\alpha_{i-1,t-1}(m) & \text{if $2\leq t\leq m$ and $t\leq i$}
\\
[2mm]
0 & \text{if $t>i$ or $t>m$}
\\
[2mm]
\frac{1}{m^{i-1}} & \text{if $t=1$ and $1\leq i$}.
\end{array}
\right.
\feq
Recall that one can define the sequence of Stirling numbers of the second kind as the solution to the recursion
\beqn
\label{compa}
S_{n,k}=kS_{n-1,k}+S_{n-1,k-1},\qquad n,k\in\nn,\,k\leq n.
\feqn
A comparison with \eqref{compa} reveals that for $t\leq m,$
\beq
P_m(N_i=t)=\frac{S_{i,t}}{m^i} \frac{m!}{(m-t)!}.
\feq
In addition,
\beq
P_m(X_{i+1}=\ell|N_i=t) =
\left\{
\begin{array}{ll}
\frac{1}{m} & \quad \text{if} \quad \ell \leq t
\\
[2mm]
\frac{m-t}{m} & \quad \ell = t +1
\\
[2mm]
0 & \quad  \mbox{otherwise}.
\end{array}
\right.
\feq
Notice that some of the boxes may remain empty at the end of the algorithm's run.
\par

In order to obtain the probability distribution of $\calf_n$, we make a simple observation that $\calf_n$ has the same distribution as that of the random variable $1\leq \calj \leq n$ defined as
\beq
\calj := \min\Big(n, \max\{j|X_1<\cdots <X_j \}\Big) = \min\Big(n, \max\{j|N_j=j \}\Big).
\feq
We will consider three cases;
\par
\begin{itemize}
\item[] {\bf Case $j=1$}: By the observation above
\beqn
P(\calf_n=j)&=&E\Big(P_M(N_2=1|N_1=1)P_m(N_1=1)\Big) = E(\frac{1}{M}) \notag \\
&=& \sum_{m=1}^{\infty} \frac{m^{n-1}}{m!\ eB_n} = \frac{B_{n-1}}{B_n}. \label{rcase1}
\feqn
\par
\item[] {\bf Case $2\leq j<n$}: We consider two possibilities such that either (i) $m<j,$ where
\beq
P_m(N_j=j, N_{j+1}=j) = 0,
\feq
or (ii) $2 \leq j\leq m$, where
\beq
P_m(N_j=j, N_{j+1}=j) &=& P_m(N_{j+1}=j | N_{j}=j) \prod_{s=1}^{j-1} P_m(N_{s+1}=s+1 | N_{s}=s) \\
&=&  \frac{j}{m} \prod_{s=1}^{j-1} \frac{m-s}{m}.
\feq
Thus, for $2\leq j < n$, our observation implies
\beqn
P(\calf_n=j) &=& \sum_{m=j}^{\infty} \frac{m^n}{m!\ e B_n} \frac{j}{m} \prod_{s=1}^{j-1} \frac{m-s}{m} = \frac{j}{eB_n} \sum_{m=j}^{\infty} \frac{m^{n-j-1}}{(m-j)!} = \frac{j\Theta_{n-j-1}(j) }{B_n},  \label{rcase2}
\feqn
where we used \eqref{Blextension} for the last equality.
\par
\item[] {\bf Case $j=n$}: We again investigate two possibilities such that either $m<n$ or $m\geq n$. Similar to the previous case, the latter is the only one contributing to the sum. Hence,
\beqn
P(\calf_n=n) = E(P_M(N_n=n)) = \sum_{m=n}^\infty \frac{1}{(m-n)!\ e B_n}  = \frac{1}{B_n}. \label{rcase3}
\feqn
\par
\end{itemize}
Finally, equations \eqref{rcase1}, \eqref{rcase2}, and \eqref{rcase3} yield \eqref{Ef_n_prob}.

{\small \begin{thebibliography}{99}

\bibitem{ABK}
M. Archibald, A. Blecher, and A. Knopfmacher, Fixed points in compositions and words, {\em J. Integer Seq.} 23 (2020), Article 20.11.1.
\filbreak

\bibitem{AT}
R. Arratia and S. Tavare, The cycle structure of random permutations, {\em Ann. Probab.} 20 (1992), 1567--1591.
\filbreak


\bibitem{Bo}
M. B\'ona, On a balanced property of derangements, {\em Electron. J. Combin.} 13 (2006), \#R102.
\filbreak

\bibitem{B}
R. A. Brualdi, {\em Introductory Combinatorics}, 5th ed., Prentice-Hall, 2010.
\filbreak

\bibitem{CMS}
N.~Cakic, T.~Mansour, and R.~Smith,
\emph{Elements protected by records in set partitions},
J. Diff. Eq. Appl. \textbf{24} (2018), 1880--1893.
\filbreak

\bibitem{C}
P. J. Cameron, {\em Combinatorics: Topics, Techniques, Algorithms}, Cambridge University Press, 1994.
\filbreak


\bibitem{maybe1}
Y.~Cai and M.~A.~Readdy,
\emph{$q$-Stirling numbers: a new view},
Adv. in Appl. Math. \textbf{86} (2017), 50--80.
\filbreak

\bibitem{comtet}
L.~Comtet,
\emph{Advanced Combinatorics. The Art of Finite and Infinite Expansions},
revised and enlarged edition, D. Reidel Publishing Co., 1974.
\filbreak


\bibitem{conflitti1}
A.~Conflitti and R.~Mamede,
\emph{Gray codes and lexicographical combinatorial generation for nonnesting and sparse nonnesting set partitions},
Theoret. Comput. Sci. \textbf{592} (2015), 87--96.
\filbreak


\bibitem{DE}
E. Deutsch and S. Elizalde, The largest and the smallest fixed points of permutations, {\em European J. Combin.} 31 (2010), 1404--1409.
\filbreak

\bibitem{DFG}
P. Diaconis, J. Fulman, and R. Guralnick, On fixed points of permutations, {\em J. Algebraic Combin.} 28 (2008), Article 189.
\filbreak


\bibitem{galvin}
D.~Galvin and A.~Pacurar,
\emph{Total non-negativity of some combinatorial matrices},
2019, preprint is available at \url{https://arxiv.org/abs/1807.08658}.
\filbreak

\bibitem{HX}
G.-N. Han and G. Xin, Permutations with extremal number of fixed points, {\em J. Combin. Theory Ser. A} 116 (2009), 449--459.
\filbreak


\bibitem{KMW}
A.~Knopfmacher, T.~Mansour, and S.~Wagner,
\emph{Records in set partitions},
Electron. J. Combin. \textbf{17} (2010), \#R109.
\filbreak


\bibitem{Mb}
T.~Mansour,
\emph{Combinatorics of Set Partitions},
CRC Press, 2013.
\filbreak


\bibitem{RAT}
T.~Mansour, R. Rastegar, and A. Roitershtein,
\emph{Height of records in partitions of a set},
Available at https://arxiv.org/abs/1908.00846.
\filbreak

\bibitem{RAT2}
T.~Mansour, R. Rastegar, and A. Roitershtein,
\emph{Horizontal visibility graph of a random restricted growth sequence},
Adv. Appl. Math. 124 (2021) 102145.
\filbreak


\bibitem{stam1}
A.~J.~Stam,
\emph{Generation of random partitions of a set by an urn model},
J. Combin. Theory Ser. A \textbf{35} (1983), 231--240.
\filbreak


\bibitem{S}
R. P. Stanley, {\em Enumerative Combinatorics}, Vol. 1, Cambridge University Press, 1986.
\filbreak

\bibitem{W}
H. S. Wilf, {\em generatingfunctionology}, A. K. Peters, 1990.
\filbreak

\end{thebibliography}
}
\end{document}